\documentclass[fleqn]{mat01}
\usepackage{times,mathtimy,amssymb,latexsym,amsbsy}
\begin{document}

\setcounter{page}{191}
\firstpage{191}

\newtheorem{theore}{Theorem}
\renewcommand\thetheore{\arabic{section}.\arabic{theore}}
\newtheorem{theor}[theore]{\bf Theorem}
\newtheorem{rem}[theore]{Remark}
\newtheorem{propo}[theore]{\rm PROPOSITION}
\newtheorem{lem}[theore]{Lemma}
\newtheorem{definit}{\rm DEFINITION}
\newtheorem{coro}[theore]{\rm COROLLARY}
\newtheorem{exampl}[theore]{Example}
\newtheorem{case}{Case}

\def\corol{\trivlist \item[\hskip \labelsep{COROLLARY.}]}
\def\noteproof{\trivlist \item[\hskip \labelsep{\it Note added in Proof.}]}

\renewcommand{\theequation}{\thesection\arabic{equation}}

\title{$\pmb{L}^{\pmb{p}}$-continuity for Calder\'on--Zygmund operator}

\markboth{Q X Yang}{$L^{p}$-continuity for Calder\'on--Zygmund operator}

\author{Q X YANG}

\address{Department of Mathematics, Wuhan University, 430072 Hubei, China\\
\noindent E-mail: yangqi99@public.wh.hb.cn}

\volume{115}

\mon{May}

\parts{2}

\pubyear{2005}

\Date{MS received 22 July 2004; revised 10 December 2004}

\begin{abstract}
Given a Calder\'{o}n--Zygmund (C--Z for short) operator $T$, which
satisfies H\"ormander condition, we prove that: if $T$ maps all
the characteristic atoms to $WL^{1}$, then $T$ is continuous from
$L^{p}$ to $L^{p}(1<p<\infty)$. So the study of strong continuity
on arbitrary function in $L^{p}$ has been changed into the study
of weak continuity on characteristic functions.
\end{abstract}

\keyword{C--Z operator; characteristic atoms; $WL^{1}$;
Hardy--Littlewood maximal operator; *-maximal operator.}

\maketitle

\section{Principal theorem}

In this paper, $L^{p}(1<p<\infty)$ continuity is obtained without
assumption on $L^{2}$ continuity, but with a continuity which is
much more weaker than the continuity from $L^{2}$ to $WL^{2} - T$
is continuous from characteristic atoms to $WL^{1}$ and no
information about its adjoint is assumed; and so, an analysis
problem is changed into a geometric problem.

Let $B(u,t)$ be a ball with center $u$ and radius $t$. A linear
operator $T$, which is continuous from $S(R^{n})$ to $S'(R^{n})$,
corresponds to a kernel distribution $K(x,y)$.

\begin{definit}$\left.\right.$\vspace{.5pc}

\noindent {\rm One calls $T$ a C--Z operator or $T\in HCZ$, if $T$
satisfies the following four conditions:

\begin{enumerate}
\renewcommand{\labelenumi}{(\roman{enumi})}
\leftskip .35pc
\item Size condition:
\begin{equation}
\hskip -1.2pc\sup\limits_{x,r}\int_{r\leq |x-y|\leq 2r}
\{|K(x,y)|+|K(y,x)|\}\hbox{d}y<\infty.
\end{equation}

\item H\"{o}rmander regular condition:
\begin{equation}
\hskip -1.2pc\sup\limits_{x,x'}\int_{|x-y|\geq 2|x-x'|}
\{|K(x,y)-K(x',y)|+|K(y,x)-K(y,x')|\}\hbox{d}y<\infty.
\end{equation}

\item T1 condition:
\begin{equation}
\hskip -1.2pc T1\in \hbox{BMO},\quad T^{*}1\in \hbox{BMO}.
\end{equation}
\newpage

\item Weak bounded condition:
\begin{gather}
\begin{split}
\hskip -1.2pc |\langle Tf,g\rangle|\leq
ct^{n}(\|f\|_{\infty}+t\|\nabla f\|_{\infty})
(\|g\|_{\infty}+t\|\nabla g\|_{\infty}),\\
\hskip -1.2pc \forall u\in R^{n}, t>0,
B(u,t),f,g\in C^{1}_{0}(B(u,t)).
\end{split}
\end{gather}
\end{enumerate}}
\end{definit}

Before, $L^{p}$ continuity is often obtained under the assumption
of $L^{2}$ continuity. Whether an operator in Definition~1.1 is
continuous on $L^{2}(R^{n})$ or not is a very difficult open
problem (see \cite{8}). In this paper, $L^{p}$ continuity is
obtained under a condition which is more weaker than $L^{2}$
continuity~--~we suppose only that the given operator is
continuous from characteristic atoms to $WL^{1}$.

In this paper, $|F|$ denotes the measure of set $F$. Let $E$ be a
cube and let $N_{E}$ be the biggest integer such that
$2^{N_{E}}|E|\leq 1$. First, we introduce some definitions about
atoms.

\begin{definit}$\left.\right.$\vspace{-.5pc}

{\rm
\begin{enumerate}
\renewcommand{\labelenumi}{(\roman{enumi})}
\leftskip .35pc
\item One calls $a(x)$ an atom on $E$ or $a(x)\in A_{1}(E)$ or
$a(x)\in A_{1}$, if $\hbox{supp}\,a(x)\subset E,
\|a(x)\|_{\infty}\leq 2^{N_{E}}$ and $\int a(x)\hbox{d}x=0$.

\item One calls $a(x)$ a characteristic atom or $a(x)\in
A^{0}_{1}(E)$ or $a(x)\in A^{0}_{1}$, if there exist two sub-cubes
$G$ and $H$ in $E$ which do not intersect each other such that
$a(x)=\frac{1}{|E|}(\chi_{F}(x)-\chi_{G}(x))$, where $F$ is a
subset of cube $H$ and $|F|=|G|$.

\item One calls $a(x)$ a special atom or $a(x)\in A(E)$ or
$a(x)\in A$, if $a(x)=\frac{1}{|E|}(\chi_{F}(x)-\chi_{G}(x))$,
where $F$ and $G$ are two sub-cubes in $E$ such that $|F\cap
G|=0$ and $|F|=|G|$.
\end{enumerate}}
\end{definit}

In fact, $a(x)\in A_{1}$ is the usual
$\infty$-atom in Hardy space, $a(x)\in A$ is defined first in
\cite{7} for Besov space $B^{0,1}_{1}$. The characteristic atom set
$A^{0}_{1}$ is composed by characteristic functions and $A\subsetneqq
A_{1}^{0} \subsetneqq A_{1}$. Now we present a
definition of Lorentz space $WL^{1}$.

\begin{definit}$\left.\right.$\vspace{.5pc}

{\rm
\noindent One calls
$f(x)\in WL^{1}(E),\ \hbox{if}\ \,
\forall \lambda>0, \lambda|\{x\hbox{:}\ |f(x)|>\lambda\}\cap
E|<\infty.$}
\end{definit}

Whatever $E=R^{n}$ or not, sometimes, one denotes
$WL^{1}(E)=WL^{1}$. It is known that $WL^{1}$ is not a Banach
space and it is only a completed metric space, because its norm
does not satisfy triangle inequality and this brings some
difficulties in the study of\break continuity.

The principal theorem in this paper is the following theorem.

\begin{theor}[\!]
Given $1<p<\infty$. If $T\in HCZ${\rm ,} then the following two
conditions are equivalent{\rm :}
\begin{align}
&T\hbox{\rm :}\ L^{p}\rightarrow L^{p}.\\[.3pc]
&T\hbox{\rm :}\ A^{0}_{1}\rightarrow WL^{1}.
\end{align}
\end{theor}

The real analysis books discuss interpolation theorem (see
\cite{1,9}). Further, several years after I have proved our
Theorem~1.1 (see \cite{12}), somebody told me that Journ\'e in 1983 proved in
\cite{6} the following theorem.
\begin{theor}[\!]
$\forall 1<p<\infty${\rm ,} if $\|T\|_{H^{1}\rightarrow
L^{1}}+\|T\|_{L^{\infty}\rightarrow {\rm BMO}}<\infty${\rm ,} then
$T\hbox{\rm :}\ L^{p}\rightarrow L^{p}$.
\end{theor}

Applying the above Theorem~1.2, the principal Theorem~1.1 can be
decomposed to the following two theorems.
\begin{theor}[\!]
If $T\in HCZ${\rm ,} then {\rm (i)} $T\hbox{\rm :}\ A_{1}\rightarrow
WL^{1}$ implies $T\hbox{\rm :}\ A_{1}\rightarrow L^{1}$ and {\rm
(ii)}~$T\hbox{\rm :}\ A^{0}_{1}\rightarrow WL^{1}$ implies $T\hbox{\rm
:}\ A^{0}_{1}\rightarrow L^{1}$.
\end{theor}

\begin{theor}[\!]
If $T\in HCZ${\rm ,} then $T\hbox{\rm :}\ A^{0}_{1}\rightarrow WL^{1}$
implies $T^{*}\hbox{\rm :}\ A_{1}\rightarrow WL^{1}$.
\end{theor}

The proof for Theorems~1.3 and 1.4 will be given in \S\S3 and 4.
Here, we apply the above three theorems to prove Theorem~1.1.

If $T\in HCZ$ and $T$ satisfies condition (1.6), by applying
Theorem~1.4, $T^{*}$ satisfies condition (1.9) below, and
hence $T^{*}$ satisfies condition (1.6). Then one applies another
time Theorem~1.4 for $T^{*}$, one gets
$T,T^{*}\hbox{\rm :}\ A_{1}\rightarrow WL^{1}$.

Further, applying (i) of Theorem~1.3, one gets
$T,T^{*}\hbox{\rm :}\ A_{1}\rightarrow L^{1}$. That is to say,
$T,T^{*}\hbox{\rm :}\ H^{1}\rightarrow L^{1}$. Since the dual space of
$H^{1}$ is BMO, one applies Theorem~1.2, $\forall
1<p<\infty,T\hbox{\rm :}\ L^{p}\rightarrow L^{p}$.

Furthermore, (1.5) implies that $T$ satisfies (1.7) below, and so
$T$ satisfies (1.6). This finishes the proof of Theorem~1.1.

\setcounter{theore}{0}
\begin{rem}
{\rm We indicate here, that so much work has been done
\cite{3,4,5,7,10,11} which deals with the continuity of C--Z operators
since the famous $T1$ theorem of David and Journ\'e \cite{2}.
Note that the following conditions are more and more weaker:
\begin{align}
&T\hbox{\rm :}\ H^{1}\rightarrow L^{1}.\\[.3pc]
&T\hbox{\rm :}\ L^{1}\rightarrow WL^{1}.\\[.3pc]
&T\hbox{\rm :}\ A_{1}\rightarrow WL^{1}.\\[.3pc]
&T\hbox{\rm :}\ A\rightarrow WL^{1}.
\end{align}
\begin{enumerate}
\renewcommand\labelenumi{(\arabic{enumi})}
\leftskip .15pc
\item A famous conjecture is: If $T\in HCZ$, then $T$ is continuous
from $L^{2}$ to $L^{2}$; Meyer \cite{8} calls it C--Z conjecture
on H\"ormander condition. According to a famous result in \cite{7},
if $T\in HCZ$, then $T$ satisfies condition (1.10); but we do not
know, under the condition that $T\in HCZ$, whether (1.6) is stronger
than (1.10) or not. Hence C--Z conjecture rests still open.

\item The conclusion in Theorem~1.3 depends only on one single side
H\"ormander\break condition.

\item As a C--Z operator $T$, it is known that $L^{2}$ continuity
implies that $T$ satisfies all the conditions from (1.5) to
(1.10). In contrast, although Journ\'e proved Theorem~1.2, in our
Theorem~1.1, we do not suppose any condition on the adjoint
operator and we suppose only a geometric condition (1.6) on the
operator itself, which is weaker than condition (1.7) or condition
(1.8). In fact, there is no characteristic atom decomposition for
$H^{1}$, hence (1.6) is weaker than (1.9). Furthermore, $WL^{1}$
is not a Banach space, so (1.9) is weaker than (1.7) or (1.8).
\end{enumerate}}
\end{rem}

There are too many constants in this paper, $C$ may be different
at each occurrence; but when a constant depends on some quantity,
this constant will be specified.

\section{Preliminaries}

In this section, we establish some results about the relations
among sets, $WL^{1}$ continuity, approximation of operators and maximal
operators. Let $M$ be Hardy--Littlewood maximal operator and
$0<\delta<1$.

First, we know that Meyer has proved the following two lemmas in
chapter 7 of \cite{8}. For arbitrary Borel set $B$ and $\forall f(x)$,
denote $\theta=\sup_{\lambda>0}\lambda|\{x\hbox{\rm :}\
|f(x)|>\lambda\}\cap B|$. Then one\break has:

\setcounter{theore}{0}
\begin{lem}
$\int_{B}|f(x)|^{\delta}\hbox{\rm d}x\leq
C(n,\delta)|B|^{1-\delta}\theta^{\delta}$.
\end{lem}

\begin{proof}
Let $E_{k}=\{x\in B,|f(x)|>2^{k}\}$ and $k_{B}$ the biggest
integer satisfying $2^{k}|B|\leq\theta$. If $k<k_{B}$, then blow
up $|E_{k}|$ to $|B|$; if $k\geq k_{B}$, then blow up $|E_{k}|$ to
$2^{-k}\theta$. Hence one has
\begin{align*}
\int_{B}|f(x)|^{\delta}\hbox{d}x &\leq
C\sum\limits^{+\infty}_{-\infty}2^{k\delta}|E_{k}| \leq
C|B|\sum\limits^{k_{B}-1}_{-\infty}2^{k\delta}+C\theta\sum\limits^{+\infty}_{k_{B}}2^{-k(1-\delta)}\\[.3pc]
&=C|B|2^{k_{B}\delta}+C\theta2^{-k_{B}(1-\delta)}\leq
C(n,\delta)|B|^{1-\delta}\theta^{\delta}.
\end{align*}\vspace{-1pc}
\end{proof}
\begin{lem}
If $E^{\delta}_{N}=\{x\hbox{\rm :}\ M|f|^{\delta}(x)>2^{N\delta}\}${\rm ,}
then $|E^{\delta}_{N}|\leq C2^{-N\delta}
\int_{E^{\delta}_{N}}|f(x)|^{\delta}\hbox{\rm d}x$.
\end{lem}

\begin{proof}
Let $g(x)=f(x)|_{E^{\delta}_{N}}$, then
\begin{equation*}
\{x\hbox{\rm :}\  M|g|^{\delta}(x)>2^{N\delta}\}= \{x\hbox{\rm :}\
M|f|^{\delta}(x)>2^{N\delta}\}.
\end{equation*}
Further, for arbitrary function $g$, one has $|\{x\hbox{\rm :}\
M|g|^{\delta}(x)>2^{N\delta}\}|\leq
C_{n}2^{-N\delta}\|g\|_{L^{\delta}}^{\delta}$. Hence
$|E^{\delta}_{N}|\leq C2^{-N\delta}
\int_{E^{\delta}_{N}}|f(x)|^{\delta}\hbox{d}x$.\hfill $\Box$
\end{proof}

For arbitrary ball or cube $B$, denote $\tilde{B}$ a ball or a
cube with the same center and double diameter. Then one can find
the following result in \cite{7}.

\begin{lem}$\left.\right.$
\begin{enumerate}
\renewcommand\labelenumi{\rm (\roman{enumi})}
\leftskip .15pc
\item $\forall T\in HCZ${\rm ,} there exists constant $C,\forall u\in
R^{n}${\rm ,} for arbitrary ball or cube $B$ with center $u${\rm
,} one has{\rm :} $\forall x\notin \tilde{B},|T\chi_{B}(x)|\leq
\frac{C|B|}{|x-u|}$.

\item  If an operator $T$ satisfies conditions {\rm (1.1)} and
{\rm (1.2),} then the following two conditions are equivalent:
\setcounter{equation}{0}
\begin{gather}
T\ \hbox{satisfies conditions (1.3) and (1.4).}\\
\|T\chi_{B}(x)\|_{L^{1}(\tilde{B})}+\|T^{*}\chi_{B}(x)\|_{L^{1}
(\tilde{B})}\leq C|B|, \forall\ \hbox{ball or cube}\ B.
\end{gather}
\end{enumerate}
\end{lem}

Now, let us consider the action of approximation operators on
atoms. Let $M_{\delta}f(x)= (M|f|^{\delta}(x))^{{1}/{\delta}}$ and
let $K_{\varepsilon}(x,y)=K(x,y)|_{|x-y|\geq\varepsilon}$ be the
kernel distribution of $T_{\varepsilon}$. Then, one has

\begin{lem}
If $T\in HCZ${\rm ,} then for arbitrary cube $E${\rm ,} we have

\begin{enumerate}
\renewcommand\labelenumi{\rm (\roman{enumi})}
\leftskip .35pc
\item $\forall a(x)\in A(E)${\rm ,} one has
\begin{equation}
|T_{\varepsilon}a(x)|\leq CM_{\delta}Ta(x)+C2^{N_{E}}.
\end{equation}

\item $\forall a(x)\in A_{1}(E)$ and $T$ satisfies condition {\rm (1.9),
(2.3)} is still true.

\item $\forall a(x)\in A_{1}^{0}(E)$ and $T$ satisfies condition
{\rm (1.6), (2.3)} is still true.
\end{enumerate}
\end{lem}

\begin{proof}
The proof of the above three conclusions in Lemma~2.4 is similar,
and so one proves only (ii). Let $B_{\varepsilon}$ be the ball $\{s\hbox{\rm :}\
|x-s| \leq\frac{\varepsilon}{2}\}$ and $\tilde{B}_{\varepsilon}$
the ball with the same center and twice the radius and let
$\tilde{B}^{c}_{\varepsilon}=R^{n}\backslash
\tilde{B}_{\varepsilon}$. Then we define $f_{2}(x)$ as follows:
(a)~$f_{2}$ is zero on $\tilde{B}_{\varepsilon}$ and
(b)~$f_{2}(s)=a(s)$ outside $\tilde{B}_{\varepsilon}$. Further, we
decompose $a-f_{2}$ into two functions $f_{1}$ and $f_{I}$:
$f_{1}$ is a function whose support is on
$\tilde{B}_{\varepsilon}$ and whose integral is zero, $f_{I}$ is a
constant function on $\tilde{B}_{\varepsilon}$ and zero outside
it. Hence $a(x)$ is decomposed into three functions and
$a=f_{1}+f_{I}+f_{2}$.

Since $|f_{2}|\leq C2^{N_{E}}$ and $|Tf_{2}(s)-Tf_{2}(x)|
=|\int_{|x-y|\geq\varepsilon}K(x,y)a(y)\hbox{d}y
-\int_{|x-y|\geq\varepsilon}$ $K(s,y)a(y)\hbox{d}y| \leq
\int_{|x-y|\geq\varepsilon}|K(x,y)-K(s,y)||a(y)|\hbox{d}y$,
according to H\"{o}rmander condition (1.2), $\forall s\in
B_{\varepsilon}$, one has: $|Tf_{2}(x)-Tf_{2}(s)|\leq C2^{N_{E}}$.
Since $Ta(s)=Tf_{1}(s)+Tf_{I}(s)+Tf_{2}(s)$, then:
$|Tf_{2}(x)|\leq |Ta(s)|+|Tf_{1}(s)|+|Tf_{I}(s)|+C2^{N_{E}}$.
\end{proof}

One makes $\delta$ order integration for $s$ on $B_{\varepsilon}$
and gets
\begin{align}
|Tf_{2}(x)|&\leq CM_{\delta}Ta(x) +C\varepsilon^{-({n}/{\delta})}
\left(\int_{B_{\varepsilon}}|Tf_{1}|^{\delta}\hbox{d}s\right)^{{1}/{\delta}}\nonumber\\[.3pc]
&\quad\,+C\varepsilon^{-({n}/{\delta})}
\left(\int_{B_{\varepsilon}}|Tf_{I}|^{\delta}\hbox{d}s\right)^{{1}/{\delta}}
+C2^{N_{E}}.
\end{align}
Applying Lemma~2.3, one gets
$\|Tf_{I}\|_{WL^{1}(\tilde{B}_{\varepsilon})}\leq
C\varepsilon^{n}2^{N_{E}}$. Since
$\varepsilon^{-n}2^{-N_{E}}f_{1}$ is an atom, by assumption,
$\|Tf_{1}\|_{WL^{1}}\leq C\varepsilon^{n}2^{N_{E}}$. Applying then
Lemma~2.1, one gets
\begin{equation}
\int_{B_{\varepsilon}}|Tf_{1}|^{\delta}\hbox{d}s +
\int_{B_{\varepsilon}}|Tf_{I}|^{\delta}\hbox{d}s \leq
C|\varepsilon|^{(1-\delta)n}\varepsilon^{n\delta}2^{\delta N_{E}}
=C\varepsilon^{n}2^{\delta N_{E}}.
\end{equation}
By (2.4) and (2.5), we have $|Tf_{2}(x)|\leq
CM_{\delta}Ta(x)+C2^{N_{E}}$.

Since $Tf_{2}(x)=T_{\varepsilon}a(x)$, one gets (2.3).

In the end of this section, we consider $*$-maximal operator of
$T$. Let $T_{*}a(x) = \sup_{\varepsilon>0}
|T_{\varepsilon}a(x)|$, then one has:

\begin{lem}
If $T\in HCZ${\rm ,} then for arbitrary cube $E${\rm ,} one
has{\rm :}

\begin{enumerate}
\renewcommand\labelenumi{\rm (\roman{enumi})}
\leftskip .15pc
\item If $T$ is satisfying condition {\rm (1.9),} then there exists
constant $C, \forall a(x)\in A_{1}(E)${\rm ,} and one has{\rm :}
\begin{equation}
|Ta(x)|\leq |T_{*}a(x)|+C2^{-N_{E}}.
\end{equation}
\item If $T$ is satisfying condition {\rm (1.6),} then there exists
constant $C, \forall a(x)\in A_{1}^{0}(E)${\rm ,} and {\rm (2.6)} is
still true.
\end{enumerate}
\end{lem}

\begin{proof}
The proof of (ii) is similar to that of (i), so one proves only
(i). If $T\in HCZ$, then applying Lemmas 2.3 and 2.4, the
operators $\{T_{\varepsilon}\}_{\varepsilon}$ satisfy the
condition (2.2) uniformly. Then there exist an operator $T_{0}$
satisfying (1.1), (1.2) and (2.2), and a subsequence
$\varepsilon_{j}$ which converges to 0 such that
$T_{\varepsilon_{j}}\rightarrow T_{0}$ in the sense of norm of
(1.1), (1.2) and (2.2). According to Lemma~2.3, $T_{0}\in HCZ$.

Further, choosing two test functions such that their supports are
disjoint, one knows that the kernel distribution of $T-T_{0}$
vanishes out of the diagonal $x=y$. Then there exists a
$L^{\infty}$ function $m(x)$ such that, for arbitrary cube $E$ and
$F$ and for arbitrary function $f(x)\in L^{\infty}(F)$, one has
\begin{equation}
\langle Tf(x),\chi_{E}(x)\rangle =\lim \langle
T_{\varepsilon_{j}}f(x),\chi_{E}(x)\rangle +\langle
m(x)f(x),\chi_{E}(x)\rangle.
\end{equation}

Finally, if $T\in HCZ$ is satisfying condition (1.9), by (2.7),
there exists constant $C$ such that, for arbitrary cube $E$ and
$\forall a(x)\in A_{1}(E)$, (2.6) is true.
\end{proof}

\section{Upgrade of regularity}

As to a C--Z operator, usually, $H^{1}\rightarrow L^{1}$
continuity is obtained by $L^{2}$ continuity; but the $L^{2}$
continuity is often established by a fixed decomposition (continuous
or discrete) of the operator. But when the regularity of $K(x,y)$
is weakened to H\"{o}rmander condition (1.2), it is difficult to
revert to the operator itself from the operators which have been
decomposed \cite{8,11,12}. Further, when we try to establish
operator's continuity under wavelet basis, the continuity from
$H^{1}$ to $L^{1}$ often needs a much stronger weak regularity
than from $L^{2}$ to $L^{2}$ \cite{11}. In this section, one
proves Theorem~1.3 through upgrading $WL^{1}$ continuity to
$L^{1}$ continuity. Since the proofs for (i) and (ii) are similar,
one proves only (i).

Now one proves (i) of Theorem~1.3 in three steps.

\begin{enumerate}
\renewcommand\labelenumi{(\arabic{enumi})}
\leftskip .15pc

\item $\forall$ atom $a(x)$ on cube $E$, $\forall
N\geq N_{E}$, the following inequality is true:
\setcounter{equation}{0}
\begin{equation}
\hskip -1.2pc 2^{N}|\{T_{*}a(x)>2^{N}\}|\leq C.
\end{equation}
Applying Lemma~2.4, one has
\begin{equation*}
\hskip -1.2pc |T_{\varepsilon}a(x)|\leq
CM_{\delta}Ta(x)+C2^{N_{E}}.
\end{equation*}
That is to say, $|T_{*}a(x)|\leq CM_{\delta}Ta(x)+C2^{N_{E}}$.

\parindent 1pc Let $E^{\delta}_{N}=\{x\hbox{\rm :}\  M|Ta|^{\delta}(x)>2^{N\delta}\}$.
According to Lemma~2.2, one has
\begin{equation*}
\hskip -1.2pc |E^{\delta}_{N}|\leq
C2^{-N\delta}\int_{E^{\delta}_{N}}|Ta(x)|^{\delta}\hbox{d}x.
\end{equation*}
Since $\|Ta(x)\| _{WL^{1}}\leq C$, applying Lemma~2.1, one gets
\begin{equation*}
\hskip -1.2pc \int_{E^{\delta}_{N}}|Ta(x)|^{\delta}\hbox{d}x\leq
C(n,\delta)|E^{\delta}_{N}|^{1-\delta}.
\end{equation*}
Combining the last two inequalities, one gets
$2^{N}|E^{\delta}_{N}|\leq C$. So (3.1) is true.

\item There exists a constant $C$, $\forall$ atom $a(x)$ on cube $E$
and $\forall N>N_{E}$, the following inequality is true:
\begin{equation}
\hskip -1.2pc |\{x\hbox{\rm :}\  T_{*}a(x)>2^{N+2}\}|\leq C2^{N_{E}-N}|\{x\hbox{\rm :}\
T_{*}a(x)>2^{N}\}|.
\end{equation}

\parindent 1pc $\forall x\in R^{n}$ and let $\|x\|$ be the maximum value of the
coordinates. Let $E_{N}=|\{x\hbox{\rm :}\  T_{*}a(x)>2^{N}\}|$. Then there
exist the biggest dyadic cubes $E_{N}^{l}$ such that
$|E^{l}_{N}\cap E^{l'}_{N}|=0$ for $l\neq l'$ and $E_{N} =
\cup_{l}E_{N}^{l}$. Let $\beta_{l}$ be the center of
$E^{l}_{N}$  and $d_{l}$ the length of $E^{l}_{N}$ and let
$\tilde{E}^{l}_{N}$ be the dyadic cube which contains $E^{l}_{N}$
with double length. Hence there exists a point $\alpha_{l}$ in
$\tilde{E}^{l}_{N}$, but $\alpha_{l}\notin E_{N}$. Hence
$|T_{*}a(\alpha_{l})|\leq 2^{N}$ and $\|\alpha_{l}-\beta_{l}\|\leq
\frac{3\sqrt{2}}{2}d_{l}$. Let $F^{l}_{N}$ be the set of $x$
belonging to $E^{l}_{N}$ such that $|T_{*}a(x)|>2^{N+2}$. In order to
prove (3.2), it is sufficient to prove
\begin{equation}
\hskip -1.2pc |F^{l}_{N}|\leq C2^{N_{E}-N}|E^{l}_{N}|,\quad
\forall N>C+N_{E}.
\end{equation}

Fix $l$, denote by $E_{N}^{l,*}$ the cube with the same center but
with 12 times diameter of $E^{l}_{N}$. Let
$a^{l}_{N}=\frac{1}{|E^{l,*}_{N}|}\int_{E^{l,*}_{N}} a(x)
\hbox{d}x$. To prove (3.3), one decomposes $a(x)$ into three
functions: $a_{p}(x)=a^{l}_{N}\chi_{E^{l,*}_{N}}(x)$,
$a_{1}(x)=(a(x)-a^{l}_{N})\chi_{E^{l,*}_{N}}(x)$ and
$a_{c}=a(x)(1-\chi_{E^{l,*}_{N}}(x))$. Let
$B_{x}=B(x,\beta_{l},\varepsilon, d_{l})= \{y\hbox{:}\ |x-y|\geq
\varepsilon\ \hbox{and}\ \|y-\beta_{l}\|\geq 6d_{l}\}$. For $u\in
R^{n}$ and $t>0$, let $B(u,t)=\{y\hbox{:}\ |y-u|<t\}$ and $B(u,t)^{c}=\{y\hbox{:}\
|y-u|\geq t\}$.

For $a_{c}(x)$, one has
\begin{align*}
\hskip -1.2pc |T_{*}a_{c}(x)-T_{*}a_{c}(\alpha_{l})|
&\leq\sup\limits_{\varepsilon>0}
|T_{\varepsilon}a_{c}(x)-T_{\varepsilon}a_{c}(a_{l})|\\[.3pc]
&=\sup\limits_{\varepsilon>0}
\left|\int_{B(x,\beta_{l},\varepsilon, d_{l})} K(x,y)a(y)\hbox{d}y\right.\\
&\quad\, \left.-\int_{B(\alpha_{l},\beta_{l},\varepsilon,
d_{l})} K(\alpha_{l},y)a(y)\hbox{d}y\right|\\[.3pc]
&=\sup\limits_{\varepsilon>0}I^{\varepsilon}.
\end{align*}
$\forall x\in E^{l}_{N}$, one considers two cases:

\bigskip
\noindent
$(\dagger)$ If $|\varepsilon|\leq
(6-\frac{3\sqrt{2}}{2})d_{l}$, then
$B_{x}=B_{\alpha_{l}}=\{y\hbox{\rm :}\  \|y-\beta_{l}\|\geq
6d_{l}\}$. Hence
$I^{\varepsilon}=|\int(K(x,y)-K(\alpha_{l},y))a_{c}(y)\hbox{d}y|\leq
\int |K(x,y)-K(\alpha_{l},y)||a_{c}(y)|\hbox{d}y$. Applying
H\"{o}rmander condition (1.2), one gets $|I^{\varepsilon}|\leq
C2^{N_{E}}$.

\noindent
$(\ddagger)$ If
$|\varepsilon|>(6-\frac{3\sqrt{2}}{2})d_{l}$, then the symmetry
difference of $B_{x}$ and $B_{\alpha_{l}}$ is $SD= B_{1} \cup
B_{2}$ where $B_{1}=B(x,\varepsilon)^{c}\cap B(\beta_{l}, 6d_{l})
\cap B(\alpha_{l},\varepsilon)$ and $B_{2}=
B(\alpha_{l},\varepsilon)^{c}\cap B(\beta_{l}, 6d_{l})$ $\cap
B(x,\varepsilon)$. Hence
\begin{align*}
\hskip -1.7pc I^{\varepsilon} &\leq \left|\int_{B_{x}}
(K(x,y)-K(\alpha_{l},y)) a(y)\hbox{d}y\right| + \int _{SD}
|K(\alpha_{l},y)||a(y)|\hbox{d}y\\[.3pc] &= I_{1}+I_{2}.
\end{align*}
As to $I_{1}$, let $C_{1}= \{y\hbox{\rm :}\  |x-y|\geq 2|x-\alpha_{l}|\}$ and
$C_{2}=\{y\hbox{\rm :}\  (\frac{3\sqrt{2}}{2}-\frac{3}{4})|x-\alpha_{l}|\leq
|x-y|\leq 2|x-\alpha_{l}|\}$. Then we have $B_{x}\subset C_{1}
\cup C_{2}$ and $I_{1}\leq \int_{C_{1}}
|K(x,y)-K(\alpha_{l},y)| |a(y)|\hbox{d}y + \int_{C_{2}}
|K(x,y)-K(\alpha_{l},y)| |a(y)|\hbox{d}y$. We apply H\"ormander
condition (1.2) to the first part and apply size condition (1.1)
to the second part, we get $I_{1}\leq C 2^{N_{E}}$. As to $I_{2}$,
we apply simply size condition (1.1), and get $I_{2}\leq C
2^{N_{E}}$. Hence one has $|I^{\varepsilon}|\leq C2^{N_{E}}$.

Combining the conclusions in the above two cases $(\dagger)$ and
$(\ddagger)$, one has
\begin{equation}
\hskip -1.2pc |T_{*}a_{c}(x)-T_{*}a_{c}(\alpha_{l})|\leq
C2^{N_{E}}.
\end{equation}

Now one proves
\begin{equation}
\hskip -1.2pc |T_{*}a_{c}(\alpha_{l})|\leq|T_{*}a(\alpha_{l})|+
C2^{N_{E}}.
\end{equation}
First, according to the definition of $T_{*}a_{c}(\alpha_{l})$,
one has
\begin{equation*}
\hskip -1.2pc T_{*}a_{c}(\alpha_{l})=\sup\limits_{\varepsilon>0}
\left|\int_{{\tiny\begin{array}{c} |\alpha_{l}-y|\geq\varepsilon\\
\|y-\beta_{l}\|\geq 12d_{l}
\end{array}}}
K(\alpha_{l},y)a(y)\hbox{d}y\right | =\sup\limits_{\varepsilon
>0}I^{\varepsilon}.
\end{equation*}

When $\varepsilon \geq d_{l}$, one has
\begin{align*}
\hskip -1.2pc I^{\varepsilon} \leq \left
|\int_{{\tiny\begin{array}{c} |\alpha_{l}-y|\geq\varepsilon\\
\|y-\beta_{l}\|\geq 12d_{l}
\end{array}}}
\!K(\alpha_{l},y)a(y)\hbox{d}y
-\int_{|\alpha_{l}-y|\geq\varepsilon}\!K(\alpha_{l},y)a(y)\hbox{d}y\right|
+T_{*}a(\alpha_{l}).
\end{align*}
Reasoning as above, one has, $I^{\varepsilon} \leq
T_{*}a(\alpha_{l})+C2^{N_{E}}$.

When $\varepsilon < d_{l}$, one has
\begin{equation*}
\hskip -1.2pc I^{\varepsilon} = \left|\int_{\|\beta_{l}-y\| \geq
8d_{l}}K(\alpha_{l},y)a(y)\hbox{d}y\right|.
\end{equation*}
Reasoning as above, one has, $I^{\varepsilon} \leq
T_{*}a(\alpha_{l})+C2^{N_{E}}$. Hence (3.5) is true.

Combining (3.4) and (3.5), one gets
\begin{equation}
\hskip -1.2pc \forall x\in E^{l}_{N},\quad |T_{*}a_{c}(x)|\leq
2^{N}+C2^{N_{E}}.
\end{equation}
According to Lemma~2.3, the $WL^{1}$ norm of $Ta_{p}(x)$ is not
greater than $C2^{N_{E}}|E^{l}_{N}|$. According to the assumption
(i) of Theorem~1.3, the $WL^{1}$ norm of $Ta_{1}(x)$ is not
greater than $C2^{N_{E}}|E^{l}_{N}|$. Applying (3.1), $\forall
N\geq N_{E}$, the following inequality is true:
\begin{equation}
\hskip -1.2pc 2^{N}|\{x\hbox{\rm :}\  T_{*}(a_{1}+a_{p})(x)>2^{N}\}|\leq
C2^{N_{E}}|E^{l}_{N}|.
\end{equation}

Let $E'_{l}=\{x\in E^{l}_{N},T_{*}(a_{1}+a_{p})(x)>2^{N+1}\}$,
then $|E'_{l}|\leq C2^{N_{E}-N}|E^{l}_{N}|$. Since
$2^{N+1}>C2^{N_{E}}|E^{l}_{N}|$, by (3.6), $F^{l}_{N}\subset
E'_{l}$, hence (3.3) is true.

\item Repeat applying (3.2), one gets
\begin{equation*}
\hskip -1.2pc \sum\limits_{N\geq C+N_{E}}2^{N}|\{x\hbox{\rm :}\
T_{*}a(x)>2^{N}\}|<\infty.
\end{equation*}
Furthermore, $\forall N\geq C+N_{E}$, and applying Lemma~2.5, one
gets
\begin{equation}
\hskip -1.2pc \{x\hbox{\rm :}\  |Ta(x)|>2^{N}\}\subset\{x\hbox{\rm :}\  T_{*}a(x)>2^{N+1}\}.
\end{equation}

Hence one gets
\begin{equation}
\hskip -1.2pc \sum\limits_{N\geq C+N_{E}}2^{N}|\{x\hbox{\rm :}\
Ta(x)>2^{N}\}|<\infty.
\end{equation}

Let $E^{d}$ be the cube which has the same center as $E$ but whose
edge has double length and $E^{c}=R^{n}\backslash E^{d}$. Applying
then H\"{o}rmander condition (1.2), $Ta(x)\in L^{1}(E^{c})$ one
gets
\begin{equation}
\hskip -1.2pc \sum\limits_{N<C+N_{E}}2^{N}|\{x\hbox{\rm :}\
Ta(x)>2^{N}\}|<\infty.
\end{equation}

By (3.9) and (3.10), $Ta(x)\in L^{1}$. Hence
$T\hbox{\rm :}\ A_{1}\rightarrow L^{1}$.
\end{enumerate}

\section{From characteristic atom to atom}

For arbitrary cube $E_{N}$, let $E^{k}_{N}$ be the cube with the
same center as $E_{N}$ but with $2^{k}$ times length. In this
section, one proves Theorem~1.4 by contradiction.

If Theorem~1.4 is not true, then $\forall N>0$, there exists cube
$E_{N}$, atom $a_{N}(x)\in A_{1}(E_{N})$ and integer $S_{N}$ such
that $2^{S_{N}}|F_{N}|\geq N$, where
$F_{N}=\{T^{*}a_{N}(x)>2^{S_{N}}\}\cap E_{N}^{1}$. By
contradiction hypothesis, \setcounter{equation}{0}
\begin{equation}
|\langle T^{*}a_{N}(x),\chi_{F_{N}}(x)\rangle|\geq N.
\end{equation}

One chooses a cube $G_{N}$ which is contained in
$E_{N}^{3}-E_{N}^{2}$ and with the same measure of $F_{N}$. Let
$b_{N}(x)=\frac{1}{|E_{N}|}(\chi_{G_{N}}(x)-\chi_{F_{N}}(x))$.
Applying then condition (1.1), one has
\begin{equation}
|\langle T^{*}a_{N}(x),\chi_{G_{N}}(x)\rangle|\leq C.
\end{equation}
Hence one has
\begin{equation}
|\langle T^{*}a_{N}(x),b_{N}(x)\rangle|\geq \frac{N-C}{|E_{N}|}.
\end{equation}

According to the assumption of Theorem~1.4, $Tb_{N}(x)\in WL^{1}$.
Applying then (ii) of Theorem~1.3, one gets that $Tb_{N}(x)\in
L^{1}$. Hence one has
\begin{equation*}
|\langle Tb_{N}(x),a_{N}(x)\rangle|\leq C\|a_{N}(x)\|_{\infty}.
\end{equation*}

Since $|\langle T^{*}a_{N}(x),b_{N}(x)\rangle|= |\langle
Tb_{N}(x),a_{N}(x)\rangle|$, then we have
\begin{equation}
|\langle T^{*} a_{N} (x), b_{N} (x) \rangle| \leq C \|a_{N} (x)
\|_{\infty}.
\end{equation}

Since $\|a_{N}(x)\|_{\infty}|E_{N}|\leq C$, then there is a
contradiction between (4.3) and (4.4).

\section*{Acknowledgements}

This research is supported by the NSFC of China (No. 10001027),
the innovation fund of Wuhan University and the subject construction
fund of Mathematics and Statistic School, Wuhan University.

\end{document}